\documentclass[12pt,reqno]{amsart}
\usepackage{amssymb,amsmath,amsthm,newlfont,enumerate}
\newtheorem{thm}{Theorem}[section] 
\newtheorem{prop}[thm]{Proposition}
\newtheorem{cor}[thm]{Corollary} 
\newtheorem{defi}[thm]{Definition}   
\newtheorem{lem}[thm]{Lemma} 
\newtheorem{rem}[thm]{Remarks}

\newcommand{\cH}{{\mathcal{H}}}

\newcommand{\cD}{{\mathcal{D}}}

\newcommand{\cA}{{\mathcal{A}}}

\newcommand{\hH}{\mathrm{H^2 }}

\newcommand{\IL}{\mathrm{L^2}}

\newcommand{\Hol}{\mathrm{Hol}}
\newcommand{\CondI}{\mathrm{(I)}}

\newcommand{\CondII}{\mathrm{(II)}}

\newcommand\TT{\mathbb{T}}
\newcommand\DD{\mathbb{D}}

\newcommand{\biindice}[3]%
{%

\begin{array}[t]{c}
{\displaystyle #1}\\
{\scriptstyle #2}\\
{\scriptstyle #3}
\end{array}

}

\begin{document}
\title[Compact composition operators]
{Compact composition operators on weighted Hilbert spaces of analytic functions}

\author[Kellay]{K. Kellay}
\address{CMI\\LATP\\Universit\'e de Provence\\
39, rue F. Joliot-Curie\\13453 Marseille cedex 13\\France}
\email{kellay@cmi.univ-mrs.fr}

\author[Lef\`evre]{P. Lef\`evre}
\address{UArtois, Laboratoire de  Math\'{e}matiques de Lens EA 2463
\\Universit\'{e} Lille Nord de France \\ F--59000 Lille \\ France }
\email{pascal.lefevre@euler.univ-artois.fr}

\thanks{The research of the first author was partially supported by 
ANR Dynop. Part of this work was done during a visit of the second named author 
at LATP}
\date{\today}

\maketitle
\begin{abstract} We characterize the  compactness  of  composition operators; in term of generalized Nevanlinna counting functions, on  a large class of Hilbert spaces of analytic functions, which can be viewed between the Bergman and the Dirichlet spaces\end{abstract}

\section{Introduction}
Let $\DD$ be the unit disk on the complex plane.  Given a positive integrable function $\omega\in C^2[0,1)$,  we extend it by $\omega(z) =
\omega (|z|)$,  $z\in \DD$, and call such $\omega$ a weight function.
We denote by  $\cH_\omega$ the space consisting of analytic functions $f$  on $\DD$ such that
$$ \|f'\|_{\omega}^{2} :=\int_\DD |f'(z)|^2\omega(z) dA(z) \,<\, \infty.$$
where $dA(z)=dxdy/\pi$ stands for the normalized area measure in $\DD$. The space $\cH_\omega$ is endowed with the norm
$$\|f\|_{\cH_\omega}^2:=|f(0)|^2+ \|f'\|_{\omega}^{2}.$$
A simple computation shows that a function $f(z)=\sum_{n=0}^\infty a_n z^n$
belongs to $\cH_\omega$
if and only if
$$
\|f\|_{\cH_\omega}^{2}=\sum_{n\geq 0}|a_n|^2 w_n< \infty,
$$
where $\omega_0=1$ and
$$\omega_n=2 n^2\int_{0}^{1}r^{2n-1}\omega(r)dr,\qquad n\geq 1.$$

{\sl Examples:} Let $\alpha>-1$, $\omega_\alpha(r)=(1-r^2)^\alpha$ and denote $\cH_{\omega_\alpha}$ by $\cH_{\alpha}$ . The Hardy space $\hH$ can be identified  with $\cH_{1}$. 
The Dirichlet space $\cD_\alpha$ is precisely $\cH_{\alpha}$ for $0\leq \alpha <1$ and $\cH_{0}$  corresponds to classical Dirichlet space $\cD$. 

Finally, the Bergman spaces  $\cA_\alpha^2(\DD):=\Hol(\DD)\cap \IL(\DD, (1-|z|^2)^\alpha dA(z))$ can be identified with  $\cH_{\alpha+2}$.
\medskip

\begin{defi}
We assume that $\omega$ is a weight function, with the following properties

\begin{description}

 \item[($W_1$)] $\omega$  is non-increasing,

 \item[($W_2$)]  $\omega(r) (1-r)^{-(1+\delta)}$ is non-decreasing for some $\delta>0$,

\item[($W_3$)] $\displaystyle \lim_{r\to 1-}\omega(r)=0$,

\item[($W_4$)] One of the two properties of convexity is fulfilled

$$\left\{\begin{array}{lll}
\displaystyle (W_{4}^{\CondI}):&& \text{  $\omega$ is convex and $\displaystyle \lim_{r\to 1}\omega'(r)=0$},\\
\text{or}&\\
\displaystyle (W_{4}^{\CondII}):&&\text{  $\omega$  is  concave. }\\
\end{array}
\right.
$$
\end{description}

Such a weight  $\omega$ is called admissible.  
\end{defi}

Sometimes, we are going to be more specific: if  $\omega$ satisfies conditions ($W_1$)--($W_3$) and ($W_4^{\CondI}$) (resp. ($W_4^{\CondII}$)), we shall say that $\omega$ is $\CondI$-admissible (resp.  $\CondII$-admissible). 

{\sl Examples:} point out that $\CondI$-admissibility corresponds to  the case $H^2\subsetneq \cH_\omega\subset \cA^2_{\alpha}(\DD)$ for some $\alpha>-1$, whereas $\CondII$-admissibility corresponds to the case  $\cD\subsetneq \cH_\omega\subseteq H^2$.

The Nevanlinna counting functions shall play a key role in our work (see \cite{LLQR1} or \cite{LLQR2} for recent results on the classical Nevanlinna counting function and the quadratic Nevanlinna counting function).

\begin{defi} Let $\varphi\in \Hol(\DD)$ such that $\varphi(\DD)\subset \DD$.

The generalized counting Nevanlinna function associated to $\omega$ is defined for every $z\in\DD\setminus\{\varphi(0)\}$ by
$$N_{\varphi,\omega}(z)=\biindice{\sum}{\varphi(a)=z}{a\in\DD}\omega(a).$$

\end{defi}

Point out that $N_{\varphi,\omega}(z)=0$ when $z\notin\varphi( \DD)$. By convention, we define $N_{\varphi,\omega}(z)=0$ when $z=\varphi(0)$. When $\omega(r)=\omega_0(r)\sim \log 1/r$, $N_{\varphi,\omega_0}=N_{\varphi}$ is the usual Nevanlinna counting function associated to $\varphi$.

In this note, we study the composition operators on $\cH_\omega$. The composition operator by $\varphi$ is defined as follows 
$$C_\varphi(f)=f\circ\varphi,\qquad f\in \cH_\omega .$$

The main result of the paper (Theorem \ref{mainthm} below) will concern compactness of $C_\varphi$. Nevertheless, before proving this result, we have to ensure the boundedness of $C_\varphi$. If  $\varphi$ is a holomorphic map on the unit disk $\DD$ into itself, it is an easy consequence of Littlewood's subordination principle (see \cite{S} for instance) that the composition operator with  $\varphi$,  induces a bounded operator $C_\varphi$ on  $\cH_\omega$ for $\CondI$--admissible weight $\omega$.

For the case of  $\CondII$--admissible weight we have

\begin{thm}\label{thmborne} Let $\omega$ be a $\CondII$--admissible weight and $\varphi\in \cH_\omega$. Then $C_\varphi$ is bounded on $\cH_\omega$
if and only if
\begin{equation}
\label{condborne}
\displaystyle \sup_{|z|< 1}\frac{N_{\varphi,\omega}(z)}{\omega(z)}<\infty
\end{equation}
\end{thm}
Note that the preceding theorem is also valid under ($W_4^{\CondI}$), but \eqref{condborne} is always fulfilled: either by a simple computation (see Remark \ref{borneI}), or using the fact that we already proved the boundedness in this case.

\vspace{1em}

The following theorem  generalizes the previously known results of \cite[Theorem 2.3, Corollary 6.11]{S1} or \cite{CS}, on Hardy and Bergman spaces see also Corollary \ref{maintcor}.  

\begin{thm}\label{mainthm} Let $\omega$ be an admissible weight and $\varphi\in \cH_\omega$. Then $C_\varphi$ is compact on $\cH_\omega$ if and only if
\begin{equation}
\label{condcompacite}
\displaystyle \lim_{|z|\to 1^-}\frac{N_{\varphi,\omega}(z)}{\omega(z)}=0
\end{equation}
\end{thm}
\medskip

Obviously, condition \eqref{condcompacite} implies the boundedness of $C_\varphi$ on $\cH_\omega$ for the admissible weight. Theorem \ref{mainthm} is the limiting case  for  the  $\cH_\alpha$ for  $\alpha>0$ and the Dirichlet space  is a degenerate case (Theorem \ref{thmborne} does not apply for the Dirichlet space).  Indeed,  Theorem \ref{mainthm} asserts that $C_\varphi$ is compact on $\cD_\alpha:=\cH_\alpha$ for $0<\alpha<1$ if and only if \eqref{condcompacite} is satisfied, i.e.
$$N_{\varphi,\alpha}(z):=\sum_{\varphi(w)=z}(1-|w|^2)^{\alpha}=o((1-|z|^2)^\alpha).$$
Note that $N_{\varphi,0}(z)$ is just the multiplicity $n_\varphi(z)$ of $\varphi$ at $z$.

Let us recall that Zorboska showed in \cite{Zo} (see also \cite{CS}) that, for $\varphi\in \cD_\alpha$ where $0\leq \alpha <1$, $C_\varphi$ is bounded  on $\cD_\alpha$ if and only if $N_{\varphi,\alpha} dA(z)$ is a Carleson measure for $\cA_\alpha(\DD)$ and $C_\varphi$ is compact on $\cD_\alpha$ if and only if $N_{\varphi,\alpha} dA(z)$ is a vanishing Carleson measure for $\cA_\alpha(\DD)$. More explicitly, for $0\leq \alpha<1$ and for all $\zeta\in \TT$,  we have 
$$\left\{\begin{array}{lll}
\text{ $C_\varphi$ is bounded on $\cD_\alpha$}&\iff\displaystyle  \sup_{\delta> 0}\frac{1}{\delta^{2+\alpha}}\int_{\{|z-\zeta|<\delta\}}N_{\varphi,\alpha}(z)dA(z)<\infty,\\
&\\
\text{ $C_\varphi$ is compact on $\cD_\alpha$}&\iff\displaystyle  \lim_{\delta\to 0}\frac{1}{\delta^{2+\alpha}}\int_{\{|z-\zeta|<\delta\}}N_{\varphi,\alpha}(z)dA(z)=0.\\
\end{array}
\right.
$$
We are going to recover these results for $\alpha>0$ as simple consequence of our results (see Theorem \ref{ThZ}).

\bigskip

There is another approach on the subject: given a continuous function $\sigma:[0,1)\to (0,\infty)$ such that    $\sigma \in \textrm{L}^1(0,1)$, we can consider the weighted Bergman space $$\cA_\sigma^2(\DD):=\Hol (\DD) \cap \IL(\DD, \sigma dA )$$ consisting  of analytic functions in $\DD$ and square area integrable with respect to the weight $\sigma$. The space  $\cA_\sigma^2(\DD)$ is  equipped with the norm $$ \|f\|_{\sigma}=\Big(\int_\DD |f(z)|^2\sigma(z)dA(z)\Big)^{1/2}.$$

If  $\varphi$ is a holomorphic map on the unit disk $\DD$ into itself, by Littlewood's subordination principle, the composition operator $C_\varphi$ is  bounded on  $\cA_\sigma^2(\DD) $.  A simple computation shows that  function $f(z)=\sum_{n=0}^\infty a_n z^n$
belongs to $\cA_\sigma^2(\DD)$
if and only if
$$
\|f\|_{ \sigma}^{2}=\sum_{n\geq 0}|a_n|^2 \sigma_n< \infty,
$$
where 
$$\sigma_n=2\int_{0}^{1}r^{2n+1}\sigma(r)dr,\qquad n\geq 0.$$
We associated to $\sigma$ the weight given by
$$\omega_\sigma(r)=\int_r^1(t-r)\sigma(t)dt.$$
Point out that $\displaystyle \lim_{r\to 1-}\omega'_\sigma(r)=0$ since $\sigma\in L^1$.

We have
$$
\frac{\sigma_{n+1}}{(1+n)^2}\asymp \int_{0}^{1}r^{2n+1}\omega_\sigma(r)dr\qquad n\geq 0.
$$
Therefore for every $f\in \cA_\sigma^2(\DD)$, we  have
$$\|f\|_\sigma^2 \asymp|f(0)|^2+\|f'\|_{\omega_\sigma}^{2}.$$
So $\cA^2_\sigma(\DD)= \cH_{\omega_\sigma}$. 

Moreover, it is worth pointing out that the weight $\omega_\sigma$ always verifies ($W_1$), ($W_3$) and ($W_4^\CondI$). So that, to know whether $\omega_\sigma$ is 
$\CondI$--admissible (and so apply the main theorem) is equivalent to know whether $\omega_\sigma$ verifies ($W_2$) or not. We have the following corollary
\begin{cor}\label{maintcor}
Let  $\varphi\in \Hol(\DD)$ such that $\varphi(\DD)\subset \DD$.  Let $\sigma$ be a weight such that  $\omega_\sigma$ is $\CondI$--admissible. 
Then $C_\varphi \text{ is compact on  } \cA_\sigma^2(\DD)$ if and only if  $$\lim_{|z|\to 1-}\frac{N_{\varphi,\omega_\sigma}(z)}{\omega_\sigma(z)}=0.$$
\end{cor}

{\sl Examples.} Note that if  $\alpha>-1$ and  $\sigma_\alpha(r)=(1-r)^{\alpha}$, then  $\omega_{\sigma_\alpha}=\omega_{\alpha+2}$.  The composition operators $C_\varphi$ is compact on
$\cA_\alpha(\DD)=\cA^{2}_{\sigma_\alpha}(\DD)$ if and only if
\begin{equation}\label{condcompaciteang}
\displaystyle \lim_{|z|\to 1^-}\frac{1-|\varphi(z)|}{1-|z|}=\infty
\end{equation}
(see \cite{CS, S1}). The condition \eqref{condcompaciteang} means that $\varphi$ does not have a finite angular derivative at any point of $\partialÊ\DD$. The compactness of $C_\varphi$ on $\hH$ implies  \eqref{condcompaciteang}, but the angular derivative condition \eqref{condcompacite} is no longer sufficient for the compactness of $C_\phi$ on $\hH$   for the general case but still sufficient for finitely valent  symbol (see \cite{CS}). Recall that $\varphi$ is finitely valent when $\displaystyle\sup_{z\in\DD}n_\varphi(z)<\infty$. We have the following corollary which involves a condition which can be viewed as a generalization of the condition \eqref{condcompaciteang}


\begin{cor}\label{corANG}
Let $\sigma$ be an admissible weight and  $\varphi\in \Hol(\DD)$ such that $\varphi(\DD)\subset \DD$. 

\begin{itemize} 
\item If $C_\varphi \text{ is compact on  } \cH_\omega$, then $\displaystyle\lim_{|z|\to1^-}\frac{\omega(z)}{\omega(\varphi(z))}=0.$

\item This is actually an equivalence, when $\varphi$ is a finitely valent holomorphic function from the disk to itself.
\end{itemize}
\end{cor}


Another example where Corollary \ref{maintcor} applies, and which (as far as we know) was not treated before in the litterature, is the following limiting case:
$$\sigma(r)=\Big((1-r^2) \log \frac{e}{1-r^2}\log\log\frac{e_2}{1-r^2}\cdots \big(\log_p\frac{e_p}{1-r^2}\big)^2 \Big)^{-1},$$
here $\log_1 x=\log x$, $\log_{k+1} x=\log\log_k x$, $e_1=e$ and $e_{k+1}=e^{e_k}$.

For this weight, it is easy to see that $\displaystyle \omega_\sigma(r)\asymp (1-r^2)\Big(\log_p \frac{e_p}{1-r^2}\Big)^{-1}$ and $\sigma_n\asymp 1/\log_p n$ so that we are closer to the Hardy space than any classical weighted Bergman space $\cA_\alpha(\DD)$, where $\alpha>-1$.

Here and in all the following, $f\asymp g$ means that there exist some constants $\alpha$, $\beta>0$ such that $\alpha f\le  g\le\beta f$.

\section{Proofs of theorems \ref{thmborne} and \ref{mainthm}}

Let $q_\lambda$ denote the automorphism of the unit disc given by
$$q_\lambda(z)=\displaystyle\frac{\lambda-z}{1-\overline{\lambda}z}, \qquad z\in \DD.$$

Consider the function $\phi=q_{\varphi(0)}\circ\varphi$.  Then $\phi:\DD\to \DD$ is analytic, $\phi(0)=0$ and $C_{q_{\varphi(0)}}$ is bounded (it suffices to make the regular change of variable). Note that  $C_\phi=C_\varphi C_{q_{\varphi(0)}}$ and since $\varphi= q_{\varphi(0)}\circ\phi$, we also have $C_\varphi=C_\phi C_{q_{\varphi(0)}}$. Therefore, $C_\varphi$ is bounded if and only if  $C_\phi$ is bounded. As well, $C_\varphi$ is compact if and only if  $C_\phi$ is compact. 

On the other hand, we have to check that the same invariance occurs on Nevanlinna counting functions, but this is an easy  consequence of the following remark: 
$$N_{q_{\varphi(0)}\circ\varphi,\omega}(z)= \sum_{q_{\varphi(0)}\circ\varphi(a)=z}\omega(a)=\sum_{\varphi(a)=q_{\varphi(0)}(z)}\omega(a)=N_{\varphi,\omega}(q_{\varphi(0)}(z))$$

At last, we can replace $\omega\big(q_{\varphi(0)}(z)\big)$ by $\omega(z)$ in the conclusion thanks to the following remark

\begin{lem}\label{poidsequiv}
If $\omega$ satisfies ($W_1$) and ($W_2$) then there exists $C>0$ such that 
$$\frac{1}{C}\omega(z)\leq \omega (q_{\varphi(0)}(z))\leq C\omega(z), \qquad z\in \DD$$
\end{lem}
\begin{proof}
Set $q_{\varphi(0)}(z)=\zeta$ and suppose that $|\zeta|\geq |z|$, By ($W_1$), we have  $\omega(\zeta)\leq \omega(z)$ and by  ($W_2$)  we get
$$\frac{\omega(z)}{\omega(\zeta)}=\frac{\omega(z)}{\omega(\zeta)}\frac{(1-|z|)^{1+\delta}}{(1-|\zeta|)^{1+\delta}} 
\frac{(1-|\zeta|)^{1+\delta}}{(1-|z|)^{1+\delta}}\leq \Big(\frac{1+|\varphi(0)|}{1-|\varphi(0)|}\Big)^{1+\delta},$$
because $1+|z|\le1+|\zeta|$ and 
$$1-|\zeta|^2=\frac{(1-|z|^2)(1-|\varphi(0)|^2)}{|1-\overline{\varphi(0)}z|^2}\cdot$$

At last, if $|\zeta|\leq |z|$, since $q_{\varphi(0)}(\zeta)=z$, it suffices to permute $z$ and $\zeta$ in the former argument. 
\end{proof}

\medskip

Hence from now until the end of the proof, we assume that $\varphi(0)=0$. In order to prove the theorems, we shall need some lemmas

\begin{lem}\label{sousharI} Let $\omega$ be a weight satisfying conditions ($W_3$) and ($W_4^\CondI$). Let  $\varphi\in \Hol(\DD)$ such that $\varphi(\DD)\subset \DD$ and $\varphi(0)=0$. 

Then the generalized Nevanlinna counting function $N_{\varphi,\omega}$ satisfies the sub--mean value property : for every $r>0$ and every $z\in\DD$
such that  $D(z,r)\subset\DD\setminus D(0,1/2)$
$$ N_{\varphi,\omega}(z)\leq  \frac{2}{r^2}\int_{D(z,r)}N_{\varphi,\omega} (\zeta)dA (\zeta)$$
\end{lem}

\begin{proof}  We set $\displaystyle \frac{d^2\omega}{dt^2}=\sigma$, so
$$\omega(t)=\int_t^1(r-t)\sigma(r)dr$$
Let $\varphi_r(z)=\varphi(rz)$, we have
\begin{eqnarray*}
N_{\varphi,\omega}(z)&=&\displaystyle\sum_{\varphi(\alpha)=z}\int_{|\alpha|}^1(r-|\alpha|)\sigma(r)dr\\
&=&\int_0^1\sum_{\varphi(\alpha)=z \atop|\alpha|\le r}(r-|\alpha|)\sigma(r)dr.
\end{eqnarray*}
Since $1/2\le|z|=|\varphi(\alpha)|\le|\alpha|\le r\le1$,
$$2(r-|\alpha|)\ge \log(r/|\alpha|)\ge r-|\alpha|.$$
So
\begin{equation}\label{subI}
2N_{\varphi,\omega}(z)\ge\int_0^1 N_{\varphi_r}(z) \sigma(r)dr \ge N_{\varphi,\omega}(z).
\end{equation}

So  by \eqref{subI}, the  sub--mean value property is inherited by the generalized Nevanlinna counting function from the same property  for the  classical Nevanlinna function (see \cite{S1} 4.6).
\end{proof}

\begin{lem}\label{sousharII} Let $\omega$ be a $\CondII$--admissible weight and let  $\varphi\in \Hol(\DD)$ such that $\varphi(\DD)\subset \DD$. Then the generalized Nevanlinna counting function $N_{\varphi,\omega}$ satisfies the sub--mean value property : for every $r>0$ and every $z\in\DD$
such that  $D(z,r)\subset\DD\setminus D(0,1/2)$
$$ N_{\varphi,\omega}(z)\leq  \frac{2}{r^2}\int_{D(z,r)}N_{\varphi,\omega}(\zeta) dA(\zeta) $$
\end{lem}

\begin{proof} By Aleman formula \cite[Lemma 2.3]{A} for $\zeta,z\in \DD$, let $\widetilde{q}_\zeta(z)=q_\zeta(-z)$ we have
\begin{equation}\label{subII}
N_{\varphi,\omega}(z)=-\frac{1}{2}\int_\DD \Delta \omega (\zeta) N_{f\circ \widetilde{q}_\zeta}(z)dA(\zeta).
\end{equation}
Note that $ \Delta \omega (\zeta) \leq 0$, since $\omega$ is decreasing and concave. We conclude as in the previous lemma.
\end{proof}

For the following lemma, we need the following well--known estimation, (see \cite[Theorem 1.7]{HKZ}).
\begin{equation}\label{integral}
\int_\DD \frac{(1-|z|^2)^c dA(z)}{|1-z\overline{\lambda}|^{2+c+d}} \asymp \frac{1}{(1-|\lambda|^2)^d}, \qquad \text{if  }d>0 \text{ , }Êc>-1,
\end{equation}

\begin{lem}\label{estimation}
Let $\omega$ be a weight satisfying ($W_1$) and ($W_2$).
$$\int_\DD\frac{\omega(z)dA(z)}{|1-\bar{\lambda}z|^{4+2\delta}}\asymp \frac{\omega(\lambda)}{(1-|\lambda|^2)^{2+2\delta}}$$
\end{lem}
\begin{proof} 
Since  $\omega$ radial and  $\omega$ is non-increasing,
$$\int_{|z|>\lambda} \frac{\omega(z)dA(z)}{|1-\bar{\lambda}z|^{4+2\delta}}\leq \omega(\lambda)\int_\DD \frac{dA(z)}{|1-\bar{\lambda}z|^{4+2\delta}}\asymp
  \frac{\omega(\lambda)}{(1-|\lambda|^2)^{2+2\delta}}$$
  The last equality follows from \eqref{integral}

On the other hand,  $\omega(r)/(1-r)^{1+\delta}$ is  non-decreasing, therefore
\begin{eqnarray*}
&&\int_{|z|<\lambda} \frac{\omega(z)}{(1-|z|^2)^{1+\delta}}\frac{(1-|z|^2)^{1+\delta}}{ |1-\bar{\lambda}z|^{4+2\delta}}dA(z)\\
& \leq & \frac{\omega(\lambda)}{(1-|\lambda|)^{1+\delta}}\int_\DD\frac{(1-|z|^2)^{1+\delta}dA(z)}{|1-\bar{\lambda}z|^{4+2\delta}}\\
&\asymp&
  \frac{\omega(\lambda)}{(1-|\lambda|^2)^{1+\delta}}\frac{1}{(1-|\lambda|^2)^{1+\delta}}
  \end{eqnarray*}
  The last equality follows from again by \eqref{integral}

The proof of the minoration is straightforward.
\end{proof}
\begin{lem}\label{fonctiontest} Let $\omega$ be a weight satisfying ($W_1$) and ($W_2$). Let $\lambda\in \DD$ and let
$$f_\lambda(z)=\frac{1}{\sqrt{\omega(\lambda)}}\frac{(1-|\lambda|^2)^{1+\delta}}{(1-\overline{\lambda} z)^{1+\delta}}$$
Then
$$\|f_\lambda\|_{\cH_\omega}\asymp 1$$
\end{lem}
\begin{proof} 

On one hand, $\displaystyle f_\lambda(0)=\frac{(1-|\lambda|^2)^{1+\delta}}{\sqrt{\omega(\lambda)}}$ is bounded by $\displaystyle  \frac{2^{1+\delta}}{\sqrt{\omega(0)}}$ thanks to ($W_2$) (actually, this even converges to $0$ when $|\lambda|\to1$).
On the other hand,
$$\|f_n'\|^2_\omega\asymp \frac{(1-|\lambda|^2)^{2(1+\delta)}}{\omega(\lambda)}\int_\DD\frac{\omega(z)}{|1-\overline{\lambda} z\big|^{4+2\delta}}dA(z)$$
The result  follows then from Lemma \ref{estimation}.
\end{proof}

\subsection*{Proof of Theorem \ref{thmborne}}
\begin{proof}
Suppose that \eqref{condborne} is satisfied. The boundedness of $C_\varphi$ follows from the change of variable formula \cite{S}:
\begin{eqnarray*}
\|C_\varphi(f)\|_{\cH_\omega}^2&=&|f(\varphi(0))|^2+\int_\DD |f'(\varphi(z))|^2|\varphi'(z)|^2\omega(z)dA(z)\\
&=&|f(0)|^2+\int_{\varphi(\DD)} |f'(z)|^2N_{\varphi,\omega}(z)dA(z)\\
&\leq &|f(0)|^2+ c\int_{\DD} |f'(z)|^2\omega(z)dA(z)\asymp \|f\|_{\cH_\omega}^2.
\end{eqnarray*}
Now assume that $C_\varphi$ is bounded on $\cH_\omega$. Let $f_\lambda$ be the test function of Lemma \ref{fonctiontest}. We have
\begin{eqnarray*}
\|C_\varphi\circ f_\lambda\|_{\cH_\omega}^2&\asymp&\frac{ (1-|\lambda|^2)^{2+2\delta}}{\omega(\lambda)}\int_{\varphi(\DD)}
\frac{N_{\varphi,\omega}(z)}{|1-\overline{\lambda} z|^{4+2\delta}}dA(z)\\
&\geq& \frac{ (1-|\lambda|^2)^{2+2\delta}}{\omega(\lambda)}\int_{D(\lambda,\frac{1-|\lambda|}{2})}\frac{N_{\varphi,\omega}(z)}{|1-\overline{\lambda} z|^{4+2\delta}}dA(z)\\
&\geq& c_1 \frac{1}{\omega(\lambda)} \frac{1}{ (1-|\lambda|^2)^{2}}\int_{D(\lambda,\frac{1-|\lambda|}{2})}N_{\varphi,\omega}(z)dA(z)\\
&\geq& c_2\frac{N_{\varphi,\omega}(\lambda)}{\omega(\lambda)}
\end{eqnarray*}
where the $c_i$'s are independent from $\lambda$ and the last inequality follows from Lemma \ref{sousharII} (when $|\lambda|$ is close enough to $1$).

We conclude that
$$\sup_{\lambda\in\DD} \frac{N_{\varphi,\omega}(\lambda)}{\omega(\lambda)}\le c'\sup_{\lambda\in\DD}\|C_\varphi\circ f_\lambda\|_{\cH_\omega}^2\le c'\|C_\varphi\|^2\sup_{\lambda\in\DD}\|f_\lambda\|_{\cH_\omega}^2$$
which is bounded by hypothesis $(C_\varphi$ bounded) and thanks to Lemma \ref{fonctiontest} .

\end{proof}
\begin{rem}\label{borneI}
\end{rem}
Still assuming that $\varphi(0)=0$, if $\omega$ is $\CondI$--admissible, \eqref{condborne} is automatically satisfied. Indeed,  the classical Littlewood's  inequality, applied to the function $r^{-1}\varphi_r$, gives that $N_{\varphi_r}(z)\le\log(r/|z|)$  and so by \eqref{subI}
\begin{eqnarray*}
N_{\varphi,\omega}(z)&\le&  \int_0^1 N_{\varphi_r}(z) \sigma(r)dr\\
& =& \int_{|z|}^1 N_{\varphi_r}(z) \sigma(r)dr\\
&\le&\int_{|z|}^1 \log(r/|z|) \sigma(r)dr\\
&\le& 2\omega(z)
\end{eqnarray*}

Up to universal constants, the same inequality is valid without assuming $\varphi(0)=0$ (see Lemma \ref{poidsequiv}).



\subsection*{Proof of Theorem \ref{mainthm}}

$\Longleftarrow$  Assume that  \eqref{condcompacite} is satisfied. Let $(f_n)_n$ be a sequence in the unit ball of $ \cH_\omega$ converging to $0$ weakly. It suffices  to show that $\|C_\varphi(f_n)\|_{\cH_\omega}\to 0$ as $n\to\infty$.  The weak convergence of $f_n$ to $0$ implies that $f_n(z)\to0$ and  $f'_n(z)\to0$ uniformly on compact subsets of $\DD$.
Let $\varepsilon>0$, there exists $\rho_\epsilon\in (1/2,1)$ such that
$$N_{\varphi,\omega}(z)\le\varepsilon\omega(z),\qquad\hbox{for } \rho_\epsilon<|z|<1.$$
By  the change of variable formula
\begin{eqnarray*}
\|C_\varphi(f_n)\|_{\cH_\omega}^{2}&\asymp &|f_n(0)|^2+\|\varphi' .( f_n'\circ\varphi)\|^2_\omega\\
&=& |f_n(0)|^2+\int_{\DD} |f'_n(\varphi(z))|^2|\varphi'(z)|^2\omega(z)dA(z)\\
&=& |f_n(0)|^2+ \int_{\varphi(\DD)} |f_n'(z)|^2 N_{\varphi,\omega}(z)dA(z)\\
&\le &|f_n(0)|^2+\int_{\rho_\varepsilon\DD} |f_n'(z)|^2N_{\varphi,\omega}(z)dA(z)+\\
&&\varepsilon\int_{\varphi(\DD)\backslash \rho_\varepsilon\DD} |f_n'(z)|^2\omega(z)dA(z)\\
&\le &|f_n(0)|^2+\int_{\rho_\varepsilon\DD} |f_n'(z)|^2N_{\varphi,\omega}(z)dA(z)+\varepsilon
\end{eqnarray*}
The conclusion easy follows since $(f'_n)$ uniformly converges to $0$ on the closed disk $\rho_\epsilon\overline{\DD}$.

\vspace{1em}

$\Longrightarrow$ Let us assume that for a $\beta>0$ and a sequence $\lambda_n\in\DD$ such that $|\lambda_n|\to1^-$ we have
$$N_{\varphi,\omega}(\lambda_n)\ge\beta.\omega(\lambda_n).$$

Let 
$$f_n(z)=\frac{1}{\sqrt{\omega(\lambda_n)}}\frac{(1-|\lambda_n|^2)^{1+\delta}}{(1-\overline{\lambda_n} z)^{1+\delta}}, \qquad z\in \DD$$

By  Lemma \ref{fonctiontest}, $(f_n)_n$ is a bounded sequence on $\cH_\omega$, converging weakly to $0$. Indeed, it is uniformly converging to $0$ on compact subsets since, by ($W_2$),   $$\frac{1}{\sqrt{\omega(\lambda_n)}}(1-|\lambda_n|^2)^{1+\delta}\le\frac{2^{1+\delta}}{\sqrt{\omega(0)}}(1-|\lambda_n|)^{(1+\delta)/2}$$

 On the other hand, by the change of variable formula and Lemmas \ref{sousharI} and \ref{sousharII} we get
\begin{eqnarray*}
\|(f_n\circ\varphi)'\|^2_\omega&\asymp&\displaystyle\frac{(1-|\lambda_n|^2)^{2+2\delta}}{\omega(\lambda_n)}
\int_\DD\frac{N_{\varphi,\omega}(z) }{\big|1-\overline{\lambda_n} z\big|^{4+2\delta}}dA(z)\\
&\ge&\displaystyle\frac{c_1}{(1-|\lambda_n|^2)^2\omega(\lambda_n)}\int_{D(\lambda_n,\frac{1-|\lambda_n|}{2})}N_{\varphi,\omega}(z)dA(z)\\
&\geq & c_1\frac{N_{\varphi,\omega}(\lambda_n)}{\omega(\lambda_n)}\\
&\geq& c_2\beta 
\end{eqnarray*}
where $c_1$ and $c_2$ are positive constant independent of $n$. Thus $C_\varphi$ cannot be compact, and this finishes the proof.
\section{Applications and complements}

First let us indicate some other special cases where Theorem \ref{mainthm} applies. Let us write $\sigma=\omega''$.

\begin{prop}

1) Condition ($W_2$) is fulfilled for every classical weight $\sigma(r)=(1-r^2)^\alpha$ where $\alpha>-1$.

2) When $\sigma$ is non decreasing, then condition ($W_2$) is fulfilled by $\omega_\sigma$ with $\delta=1$.
\end{prop}

\begin{proof}
1. Take $\delta=\alpha+1$.

2. Let compute the derivative of $\displaystyle H(r)=\frac{\omega(r)}{(1-r)^{2}}\cdot$

$$H'(r)= \displaystyle\frac{2}{(1-r)^3}\int_r^1(x-c)\sigma(x)dx$$
where $c=(r+1)/2$.
So
$$H'(r)= \displaystyle\frac{2}{(1-r)^3}\int_c^1(t-c)(\sigma(t)-\sigma(2c-t))dt\ge0$$
since $\sigma$ is non decreasing.
\end{proof}

It would be interesting to compare the condition obtained for one specific weight $\omega$ with the one for another weight. In particular, it is known that the compactness on the Hardy space $H^2$ implies the compactness on classical weighted Bergman spaces. We are going to extend this result. On the other hand, it would be interesting to know when the non angular derivative condition \eqref{condcompaciteang} is still equivalent to the compactness on weighted Bergman spaces. We also have a partial result in this direction.

Before proving these results, we need some simple observations.

We associate to $\sigma$ the weight $\omega_\sigma$ and we introduce the function 
$$G(r)=\frac{\omega_\sigma(r)}{(1-r)}, \qquad r\in [0,1[$$ 
which is $C^1$ on $[0,1[$ and extends continuously at $1$ by 
$$G(1)=-\displaystyle\lim_{r\rightarrow1}\omega_\sigma'(r)=0.$$
 Moreover, when $\omega_\sigma$ verifies ($W_3$) and $\displaystyle \lim_{r\to 1}\omega_\sigma'(r)=0$, we can write
 
$$G(r)=\int_r^1\frac{\rho-r}{(1-r)}\sigma(\rho)d\rho, $$
so 
 $$G'(r)=\int_r^1\frac{\rho-1}{(1-r)^2}\sigma(\rho)d\rho.$$
Hence $G$ is a non-increasing function. Moreover: 
\begin{eqnarray*}
G''(r)&=&\frac{\sigma(r)}{(1-r)}-\int_r^1\frac{2(1-\rho)}{(1-r)^3}\sigma(\rho)d\rho\\
&=&\frac{2}{(1-r)^3}\int_r^1(1-\rho)(\sigma(r)-\sigma(\rho))d\rho
\end{eqnarray*}
Therefore $G$ is a convex function when $\sigma$ is a non increasing function.

We introduce the following condition. Here, it will be more convenient to write $\tilde\omega_\sigma(x)=\omega_\sigma(1-x)$ and $\tilde G(x)=G(1-x)$. 

We say that the weight $\omega_\sigma$ verifies condition ($\kappa$) if 

$$\displaystyle\lim_{\eta\to 0^+}\limsup_{x\to 0^+}\frac{\tilde\omega_\sigma(\eta x)}{\eta\tilde\omega_\sigma(x)}=0.$$

This is clearly equivalent to 

$$\displaystyle\lim_{\eta\to 0^+}\limsup_{x\to 0^+}\frac{\tilde G(\eta x)}{\tilde G(x)}=0$$

Observe that when $\sigma$ is a non increasing function, then Condition ($\kappa$) is fulfilled for $\omega_\sigma$: $\tilde G(0)=0$, so by convexity:  $\tilde G(\eta x)\le\eta\tilde G(x)$.
\bigskip

\begin{thm}\label{Th1}
Let  $\varphi\in \Hol(\DD)$ such that $\varphi(\DD)\subset \DD$. Let $\sigma$ be a weight such that $\omega_\sigma$ is $\CondI$-admissible.  

\begin{enumerate}[i)]
\item  The compactness of $C_\varphi$ on the Hardy spaces $H^2$ always implies the compactness on $ \cA_\sigma^2(\DD)$.

\item The compactness of $C_\varphi$ on the weighted Bergman $\cA_\sigma^2(\DD)$ always implies the compactness on the classical Bergman space $ \cA_0^2(\DD)$ (hence condition \eqref{condcompaciteang} is fulfilled).
\end{enumerate}
\end{thm}

\begin{proof}
As already explained, it suffices to treat the case $\varphi(0)=0$.

$i)$ The Schwarz lemma implies that for every $a\in\DD$ with  $\varphi(a)=z$, we have $|z|\le|a|$, hence $G(|z|)\ge G(|a|)$ where $G(r)=\omega_\sigma(r)(1-r)^{-1}$. Then 

\begin{eqnarray*}
N_{\varphi,\omega}(z)&=&\biindice{\sum}{\varphi(a)=z}{a\in\DD}\omega_\sigma(a)\\
&=&\biindice{\sum}{\varphi(a)=z}{a\in\DD}G(|a|)(1-|a|)\\
&\le& G(|z|)\biindice{\sum}{\varphi(a)=z}{a\in\DD}(1-|a|)
\end{eqnarray*}
By hypothesis $\displaystyle\biindice{\sum}{\varphi(a)=z}{a\in\DD}(1-|a|)=o(1-|z|)$ and we conclude that 
$$N_{\varphi,\omega_\sigma}(z)=o(\omega_\sigma(z)).$$ 

$ii)$ Using the function $H(r)={\omega_\sigma(r)}(1-r)^{-(1+\delta)}$ instead of $G$, the same trick works to show that when $C_\varphi$ is compact on  $\cA_\sigma^2(\DD)$, then $C_\varphi$ is compact on  $\cA_{\sigma_{1+\delta}}^2(\DD)$. But it is known to be equivalent to the compactness on the standard Bergman space $\cA ^2_0(\DD)$.
\end{proof}
\medskip

Applying the same ideas, we are able to produce a simple sufficient  test-condition to ensure that a composition operator on a weighted Bergman space is compact. We shall see below that the converse of the first assertion is false in general, for any $\CondI$-admissible weight.

\begin{thm}\label{Th2}
Let  $\varphi\in \Hol(\DD)$ such that $\varphi(\DD)\subset \DD$. Let $\sigma$ be a weight such that $\omega_\sigma$ is $\CondI$-admissible.  

\begin{enumerate}[i)]
\item  If $\displaystyle\lim_{|z|\to1^-}\frac{G(z)}{G(\varphi(z))}=0$ then $C_\varphi$ is compact on $\cA_\sigma^2(\DD)$.

\item When $\omega_\sigma$ satisfies the condition ($\kappa$), condition \eqref{condcompaciteang} implies that $C_\varphi$ is compact on $\cA_\sigma^2(\DD)$.
\end{enumerate}
\end{thm}

Hence by Theorems \ref{Th1} and \ref{Th2}, we have

\begin{eqnarray*}
& &\displaystyle\frac{\omega(z)}{\omega(\varphi(z))}=  o\Big(\frac{1-|z|}{1-| \varphi(z)|}\Big)\\ 
&\Longrightarrow&\quad C_\varphi \text{ is compact on  }\cA_\sigma^2(\DD)\\
&\Longrightarrow&\displaystyle\frac{\omega(z)}{\omega( \varphi(z))}=o(1)
\end{eqnarray*}

when $|z|$ tends to $1$.
\medskip

Of course, in the very special case of classical weighted Bergman spaces, we recover the well-known equivalence with condition \eqref{condcompaciteang}.

\begin{proof}
As already explained, it suffices to treat the case $\varphi(0)=0$.

$i)$ Let us fix some $1>\epsilon>0$. There exists some $\rho\in(0,1/2)$ such that for every $|a|>1-\rho$, we have $\displaystyle \frac{G(a)}{G(\varphi(a))}\le\varepsilon$.

 The Schwarz lemma implies that for every $a\in\DD$ with  $\varphi(a)=z$, we have $|z|\le|a|$, so that  $|a|>1-\rho$ as soon as $|z|>1-\rho$. We have then
\begin{eqnarray*}
N_{\varphi,\omega_\sigma}(z)&=&\biindice{\sum}{\varphi(a)=z}{a\in\DD} G(|a|)(1-|a|)\\
&\le&\varepsilon\biindice{\sum}{\varphi(a)=z}{a\in\DD}G\big(\varphi(a)\big)(1-|a|)\\
&=&\varepsilon G(|z|)\biindice{\sum}{\varphi(a)=z}{a\in\DD}(1-|a|)
\end{eqnarray*}
But 
$$\displaystyle\biindice{\sum}{\varphi(a)=z}{a\in\DD}(1-|a|)\le2(1-|z|)$$ and we obtain that 
 $N_{\varphi,\omega_\sigma}(z)\le2\varepsilon\omega_{\sigma}(z)$

$(ii)$ Obvious with the previous result.  
\end{proof}

We already said that the first implication in Theorem \ref{Th1}, is not an equivalence. Indeed:

\begin{cor}
Let $\sigma$ be a weight such that $\omega_\sigma$ is $\CondI$-admissible. There exists an analytic function $\varphi:\DD\to\DD$ such that

\begin{itemize} 
\item  $C_\varphi$ is compact on  $\cA_\sigma^2(\DD)$.

\item $C_\varphi$ is not compact on the classical Hardy space $H^2$.
\end{itemize}
\end{cor}
\begin{proof}
It suffices to apply both the preceding theorem and Theorem 3.1. \cite{LLQR3} with 
$$\Delta(t)=\tilde G^{-1}\Big(C\sqrt{\tilde G(t)}\Big),$$
where the numerical constant $C$ is only fixed by $\Delta(1)=1/2$; point out that $\Delta$ is non-decreasing and that $\lim_{t\to0}\Delta(t)=0$. This provides us with a Blaschke product $\varphi$ such that $\varphi(0)=0$ and 
$$1-|\varphi(z)|\ge\Delta(1-|z|)$$
Since $\varphi$ is inner, $C_\varphi$ cannot be compact on $H^2$.

On the other hand, 
$$G(|\varphi(z)|)=\tilde G(1-|\varphi(z)|)\ge\tilde G(\Delta(1-|z|))=C\sqrt{G(|z|)}.$$
Hence

$$\displaystyle\lim_{|z|\to1^-}\frac{G(z)}{G(\varphi(z))}\le\displaystyle\lim_{|z|\to1^-} \sqrt{G(|z|)}=0$$

\end{proof}

As announced in the introduction, we can recover the characterization due to Zobroska of compactness for classical weighted Dirichlet spaces. Actually we have such characterizations for every $\CondII$-admissible weight. We use the notations already introduced in the first section.

\begin{thm}\label{ThZ} Let $\omega$ be  a $\CondII$--admissible weight and let $\varphi\in \cH_\omega$. Then  
\begin{enumerate}
\item[(i)] $C_\varphi$ is bounded on $\cH_\omega$ if and only if
$$ \sup_{\delta> 0}\sup_{\zeta\in \TT}\frac{1}{\delta^{2}\omega(1-\delta)}\int_{\{|z-\zeta|<\delta\}}N_{\varphi,\omega}(z)dA(z)<\infty.$$
\item[(ii)] $C_\varphi$ is compact on $\cH_\omega$ if and only if  
$$ \lim_{\delta\to 0}\sup_{\zeta\in \TT}\frac{1}{\delta^{2}\omega(1-\delta)}\int_{\{|z-\zeta|<\delta\}}N_{\varphi,\omega}(z)dA(z)=0.$$
\end{enumerate}
\end{thm}

\begin{proof} We only proof (ii) since the proof of (i) is similar. If we assume that $C_\varphi$ is compact. The characterization (2) easily implies that 
$$ \lim_{\delta\to 0}\frac{1}{\delta^{2}\omega(1-\delta)}\int_{\{|z-\zeta|<\delta\}}N_{\varphi,\omega}(z)dA(z)=0$$
uniformly in $\zeta\in \TT$.

On the converse, let $\theta\in\DD$ (with say $|\theta|>1/2$). Let $\zeta=\theta/|\theta|\in\TT$ and $\delta>0$ such that $\theta$ is the midpoint of $[(1-\delta)\zeta,\zeta]$. In other words: $\delta=2(1-|\theta|)\in(0,1)$. Then by Lemma 2.2 
\begin{eqnarray*}\frac{1}{\delta^{2}\omega(1-\delta)}\int_{\{|z-\zeta|<\delta\}}N_{\varphi,\omega}(z)dA(z)&\geq&\frac{1}{\delta^{2}\omega(1-\delta)} \int_{D(\theta,\frac{\delta}{2})}N_{\varphi,\omega}(z)dA(z)\\
&\geq& c\frac{N_{\varphi,\omega}(\theta)}{\omega(\theta)}\cdot
\end{eqnarray*}
Letting $|\theta|$ tend to $1$,  we get the characterization (2).
\end{proof}

\end{document}